\documentstyle[12pt]{amsart}
\newtheorem{Theorem}{Theorem}[section] 
\newtheorem{Lemma}[Theorem]{Lemma}
\newtheorem{Corollary}[Theorem]{Corollary}
\newtheorem{Proposition}[Theorem]{Proposition}
\newtheorem{Example}[Theorem]{Example}

\def\reg{\operatorname{reg}} 
\def\ini{\operatorname{in}} 
\def\Gin{\operatorname{Gin}} 
\def\rk{\operatorname{rank}}
\def\core{\operatorname{core}}

\def\sk{\smallskip}
\def\aa{{\frak a}} \def\mm{{\frak m}} 
\def\qq{{\frak q}} 
\def\QQ{{\cal Q}}  

\begin{document}

\title{Constructive characterization \\ of the reduction numbers}
\author{Ng\^o Vi\^et Trung}
\address{Institute of Mathematics,  Box 631, B\`o H\^o, 10000 Hanoi, Vietnam}
\email{nvtrung@@thevinh.ncst.ac.vn}
\thanks{The author is partially supported by the National Basic Research Program of Vietnam} 

\begin{abstract} We present a constructive description of minimal reductions with a given reduction number. This description has interesting consequences on the minimal reduction number, the big reduction number, and the core of an ideal. In particular, it helps solve a conjecture of Vasconcelos on the relationship between reduction numbers and initial ideals. \end{abstract} 
\maketitle

\centerline{\small \it Dedicated to J\"urgen Herzog on the occasion of his sixtieth birthday} \vskip 0.7cm

\section*{Introduction} \sk

Let $A$ be a standard graded algebra over an infinite field $k$ with $d = \dim A$. An ideal $Q = (t_1,\ldots,t_d)$, where $t_1,\ldots,t_d$ are linear forms of $A$, is called a {\it minimal reduction} of $A$ if $k[t_1,\ldots,t_d] \hookrightarrow A$ is a Noether normalization of $A$. The reduction number of $A$ with respect to $Q$, written as $r_Q(A)$, is the maximum degree of the generators of $A$ as a graded $k[t_1,\ldots,t_d]$-module. The {\it reduction number} of $A$ is defined as the invariant
$$ r(A) := \min\{r_Q(A)|\ \text{$Q$ is a minimal reduction of
$A$}\}.$$
The above notions were originally introduced for ideals of local rings by Northcott and Rees [NR] (see also Lipman [L]). However, the local case can be passed to the graded case via the fiber ring of the given ideal. 
\sk

The reduction number has been extremely useful in the study of blowup algebras (see e.g.~[S], [A], [AHT], [GNN], [HH1], [HH2], [JK], [JU], [U]). In general, $r(A)$ can be used as a measure for the complexity of $A$. It can be compared with other important invariants of $A$ such that the degree, the arithmetic degree and the Castelnuovo-Mumford regularity (see [T1], [V1], [V2]). Calculating $r(A)$ is usually hard since we do not know which minimal reduction has this reduction number. That is a reason for why the reduction number is not well-understood. For instance, if $A = R/I$, where $R$ is a polynomial ring over $k$ and $I$ is a homogeneous ideal of $R$, Vasconcelos conjectured that $r(R/I) \le r(R/\ini(I)),$ where $\ini(I)$ denotes the initial ideal of $I$ with respect to an arbitrary term order of $R$ [V2, Conjecture 7.2]. This conjecture has been solved for generic initial ideals by Bresinsky and Hoa [BH]. \sk

In this paper we will solve the above problems. First, in Section 1 we show that $r(A)$ is the reduction number of any generic minimal reduction and that $r(A)$ can be computed by evaluating the defining ideal of $A$. Our approach is based on the comparison of the ranks of certain ``generic" matrices with the Hilbert function of the graded algebra. It resembles, in certain sense, Eakin and Sathaye's estimation for the reduction number in [ES]. In Section 2 we describe the parameter space of the minimal reductions of $A$ with a given reduction number. This description allows us to give an explicit characterization of the big reduction number
$$br(A) := \max\{r_Q(A)|\ \text{$Q$ is a minimal reduction of $A$}\},$$
which has been studied recently by Vasconcelos [V3]. Using the big reduction number we can describe the parameter space of all minimal reductions or, equivalently, all standard Noether normalizations of $A$. Moreover, we can show that the range of the reduction numbers of the minimal reductions of $A$ needs not to be a consecutive sequence of integers. \sk

In Section 3 we will settle the above conjecture of Vasconcelos in the affirmative. 
The solution is based on the construction of a flat family connecting a given ideal $I$ to its initial ideal $\ini(I)$ via an integral weight function (see Eisenbud [E, 15.8]). More generally, we prove that $r(R/I) \le r(R/\ini_{\lambda}(I)),$ where $\ini_{\lambda}(I)$ denotes the initial ideal of $I$ with respect to an arbitrary weight function $\lambda$. A crucial point in the proof is the fact that the reduction number is preserved by any transcendental extension of the base field, which follows from the characterization of $r(A)$ in Section 1. We would like to point out that Conca [C] has independently solved Vaconcelos conjecture by a completely different method and that in general one can not compare $br(R/I)$ and $br(R/\ini(I))$. \sk

In the final Section 4 we will apply our approach to study the reduction numbers of ideals in local rings with infinite residue field. Similarly as in the graded case, the reduction number is always attained by a generic minimal reduction  and the parameter space of the minimal reductions can be described explicitly. This description allows us to specify an ideal which is contained in the core (the intersection of all minimal reductions) of the given ideal. 
The core of an ideal was first studied by Rees and Sally [RS] and then by Huneke and Swanson [HS], Corso, Pollini and Ulrich [CPU], partly due to its relationship to the theorem of Brian\c con and Skoda. As noted in [HS], it is an important question to understand how the core is related to the given ideal. The main result of [RS] states that for $\mm$-primary ideals of a Cohen-Macaulay local ring with maximal ideal $\mm$, the contraction of the generic minimal reduction is contained in the core. It is shown recently in [CPU] that this containment is actually an equality and that this equality holds under certain general assumptions. Our result is a generalization of the main result of [RS] to arbitrary ideals of an arbitrary local ring. We shall see that the specified ideal coincides with the contraction of the generic minimal reduction in most of the cases considered in [CPU], and we will give an example where the core is equal to  the specified ideal but different from the contraction of the generic minimal reduction.

\section{Generic minimal reduction} \sk

Let $A = \oplus_{n\ge 0}A_n$ be a standard graded algebra over an infinite field $k$ with $d = \dim A$. The reduction number is usually defined as follows. A {\it reduction} of $A$ is a graded ideal $Q$ generated by linear forms such that $Q_n = A_n$ for all large $n$. The least non-negative integer $n$ for which $Q_{n+1} = A_{n+1}$ is called the {\it reduction number} of $A$ with respect to $Q$. It will be denoted by $r_Q(A)$. A reduction $Q$ of $A$ is {\it minimal} if $Q$ does not contain any other reduction of $A$. It is known that a reduction of $A$ is minimal if it is generated by $d$ linear forms [NR] and that it is exactly the one defined in the introduction [V1]. We would like to point out that different minimal reductions of $A$ may have different reduction numbers [BH, Example 7] (cf. [Hu, Example 3.1] for the local case). \sk

Let $A = R/I$, where $R = k[x_1,\ldots,x_m]$ and $I$ is a homogeneous ideal of $R$. Let $Q$ be an arbitrary ideal generated by $d$ linear forms of $A$. We can find a family $\alpha =  (\alpha_{ij}) \in k^{md}$ such that if we set
$$y_i = \sum_{j=1}^m\alpha_{ij}x_j\    (i = 1,\ldots,d)$$
then  $Q = (y_1,\ldots,y_d,I)/I$. We say that $Q$ is {\it parameterized} by $\alpha$. \sk

Assume that $I = (f_1,\ldots,f_v)$, where $f_j$ is a form of degree $d_j$ of $R$, $j = 1,\ldots,v$. For every integer $n \ge 0$, the vector space $(y_1,\ldots,y_d,I)_n$ is spanned by the elements $y_ig$, where $g$ is a monomial of degree $n-1$, and $f_jh$, where $h$ is a monomial of degree $n-d_j$. Let $M_n(\alpha)$ denote the matrix of the coefficients of these elements written as linear combinations of the monomials of degree $n$. Note that we always have 
$$\rk M_n(\alpha) \le {n+m-1 \choose m-1}.$$

We can easily test when $Q$ is a minimal reduction of $A$ and compute $r_Q(A)$ by means of $M_n(\alpha)$.\sk

\begin{Lemma} \label{test}  The ideal $Q$ is a minimal reduction of $A$ if and only if there is an integer $n \ge 0$ such that $\rk M_n(\alpha) = {n+m-1 \choose m-1}$. In this case,
$$r_Q(A) = \min\{n|\ \rk M_{n+1}(\alpha) = {n+m \choose m-1}\}.$$ \end{Lemma}

\begin{pf} Since $\dim (y_1,\ldots,y_d,I)_n = \rk M_n(\alpha)$, we have
$(y_1,\ldots,y_d,I)_n = R_n$ (i.e. $Q_n = A_n$) if and only if $$\rk M_n(\alpha) = \dim R_n = {n+m-1 \choose m-1}.$$
Hence the conclusion follows from the above definition of $r_Q(A)$. \end{pf} \sk
 
We may replace $R$ by the polynomial ring $R_u := k(u)[x_1,\ldots,x_m]$, where $u = \{u_{ij}|\ i = 1,\ldots,d,\ j = 1,\ldots,m\}$ is a family of $md$ indeterminates. In $R_u$ we consider the linear forms
$$z_i = \sum_{j=1}^mu_{ij}x_j\ (i = 1,\ldots,d).$$
Put  $A_u = R_u/IR_u$ and $Q_u = (z_1,\ldots,z_d,I)/IR_u$. It is well known that $Q_u$ is a minimal reduction of $A_u$. Following the terminology of [RS] we call $Q_u$ a {\it generic minimal reduction} of $A$. \sk

For a generic minimal reduction  $Q_u$ and every integer $n \ge 0$ we can define a matrix $M_n(u)$ similarly as the matrix $M_n(\alpha)$ for the ideal $Q$. The entries of $M_n(u)$ are homogeneous polynomials in $k[u]$. Note that every matrix $M_n(\alpha)$ can be obtained from $M_n(u)$ by the substitution $u = \alpha$. \sk

The reduction number $r(A)$ can be described in terms of the matrices $M_n(u)$ as follows.\sk

\begin{Theorem} \label{generic} Let $Q_u$ be a generic minimal reduction of $A$. Then
$$r(A) = r_{Q_u}(A_u) = \min\{n|\ \rk M_{n+1}(u) = {n+m \choose m-1}\}.$$ \end{Theorem}

\begin{pf} By Lemma \ref{test} we have
$$r_{Q_u}(A_u) = \min\{n|\ \rk M_{n+1}(u) = {n+m \choose m-1}\}.$$
Put $r =r_{Q_u}(A_u)$. Then $\rk M_{r+1}(u) = {r+m \choose m-1}.$ Since $k$ is an infinite field, we can find a family $\alpha$ such that 
$$\rk M_{r+1}(\alpha) = \rk M_{r+1}(u) = {r+m \choose m-1}.$$
Hence $Q$ is a minimal reduction of $A$ with $r_Q(A) \le r$. So we obtain $r(A) \le r = r_{Q_u}(A_u)$. \par
Conversely, let $s = r(A)$. Then there exists a minimal reduction $Q$ of $A$ with $r_Q(A) = s$. Applying Lemma \ref{test} we get 
$$\rk M_{s+1}(\alpha) = {s+m \choose m-1} = \dim_{k(u)}S_{s+1}.$$
By specialization, $\dim_{k(u)} (z_1,\ldots,z_d,I)_{s+1} = \rk M_{s+1}(u) \ge \dim_{k(u)}S_{s+1}.$ Therefore, $(z_1,\ldots,z_d,I)_{s+1} = S_{s+1}$ or, equivalently, $(Q_u)_{s+1} = (A_u)_{s+1}$. This implies $r_{Q_u}(A_u) \le s = r(A)$. So we can conclude that $r(A) = r_{Q_u}(A_u)$. \end{pf} 

As a consequence of Theorem \ref{generic} we can use certain evaluations of the defining ideal $I$ of $A$ to compute $r(A)$.\sk

\begin{Corollary} Let  $v = \{v_{ij}|\ i = m-d+1,\ldots,m,\ j = 1,\ldots,m-d\}$ be a set of new indeterminates and $R' = k(v)[x_1,\ldots,x_{m-d}]$. Let $I' \subset R'$ denote the ideal generated by the elements obtained from $I$ by the substitution 
$$x_i = v_{i1}x_1 + \cdots + v_{im-d}x_{m-d},\; i = m-d+1,\ldots,m.$$
Then $r(A)$ is the largest integer $n$ such that $(R'/I')_n \neq 0$. \end{Corollary} 

\begin{pf} Solving the system of linear equations 
$$\sum_{j=m-d+1}^mu_{ij}x_j = z_i - \sum_{j=1}^{m-d}u_{ij}x_j\ (i = 1,\ldots,d)$$ 
in the variables $x_{m-d+1},\ldots,x_m$ we get 
$$x_i = \sum_{j=1}^{m-d}w_{ij}x_j + \sum_{j=1}^d w_{im-d+j}z_j\ (i = m-d+1,\ldots,m),$$
where $w = \{w_{ij}|\ i = m-d+1,\ldots,m,\ j = 1,\ldots,m\}$ is a set of elements of $k(u)$. It is easy to see that $k(u) = k(w)$ so that $w$ is a set of algebraically independent elements over $k$. Put $v_{ij} = w_{im-d+j}$ for $i = m-d+1,\ldots,m$, $j = 1,\ldots,m-d$. Then $k(u) = k(v \cup v')$, where $v' = \{v_{ij}|\ i,j = m-d+1,\ldots,m\}$. Therefore,
$$S/(z_1,\ldots,z_d,I)S \cong (R'/I') \otimes_kk(v').$$
By Theorem \ref{generic}, $r(A)$ is the largest integer $n$ such that $S_n \neq (z_1,\ldots,z_d,I)_n$. Hence $r(A)$ is the largest integer $n$ such that $(R'/I')_n \neq 0$.
\end{pf} 

Another interesting consequence of Theorem \ref{generic} is the fact that the reduction number is preserved by transcendental extensions of the base field.\sk

\begin{Corollary} \label{extension} Let $A^* = A \otimes_k k(t)$, 
where $t$ is an indeterminate. Then $$r(A^*) = r(A).$$
\end{Corollary}

\begin{pf} Let $Q_u = (z_1,\ldots,z_d)$ be a generic minimal reduction of $A$. Put $A_u^* = A_u\otimes_k k(t)$ and $Q^*_u = (z_1,\ldots,z_d)A_u^*$. Then $Q_u^*$ is also a generic minimal reduction of $A_u$. By Theorem \ref{generic}, $r(A) = r_{Q_u}(A_u)$ and $ r(A^*) = r_{Q_u^*}(A_u^*) $. It is clear that $(Q_u)_n = (A_u)_n$ if and only if $(Q_u^*)_n = (A_u^*)_n$. From this it follows that $r_{Q_u}(A_u) = r_{Q_u^*}(A_u^*)$. Hence $r(A^*) = r(A).$ \end{pf}  

\section{Parameter spaces of minimal reductions}

Let $A$ be standard graded algebra over an infinite field $k$ with $d = \dim A$. As in Section 1, let $M_n(u)$ ($n \ge 0$) be the matrices associated with a generic minimal reduction of $A$. \sk

For every integer $n \ge 0$ let $J_n$ denote the homogeneous ideal of $k[u]$ generated by the ${n+m-1 \choose m-1} \times {n+m-1 \choose m-1}$ minors of the matrix $M_n(u)$ ($m$ depends on the definition of $M_n(u)$). Note that $J_n  = 0$ if one of the sizes of $M_n(u)$ is less than ${n+m-1 \choose m-1}$. Let ${\cal V}_n \subseteq {\Bbb P}_k^{md-1}$ denote the projective variety defined by the zero locus of $J_n$. \sk

We will use the projective varieties $\{{\cal V}_n\}$ to describe the parameter space of minimal reductions of $A$ with a given reduction number. \sk

\begin{Theorem} \label{given} Let $n \ge 0$ be any given integer. Let $Q$ be an ideal generated by $d$ linear forms of $A$ which is parameterized by and parameterized by $\alpha \in {\Bbb P}_k^{md-1}$. Then\par
\noindent {\rm (i)} ${\cal V}_n \supseteq {\cal V}_{n+1}$, \par
\noindent {\rm (ii)}  $Q$ is a minimal reduction of $A$ with $r_Q(A) = n$ if and only if $\alpha \in {\cal V}_n \setminus {\cal V}_{n+1}$. \end{Theorem}

\begin{pf} It is clear that $\rk M_n(\alpha) = {n+m-1 \choose m-1}$ if and only if $\alpha \not\in {\cal V}_n$. As we have already seen in the proof of Lemma \ref{test}, we can express the condition $\rk M_n(\alpha) = {n+m-1 \choose m-1}$ as $Q_n = A_n$.  Since $Q_n = A_n$ implies $Q_{n+1} = A_{n+1}$, we must have ${\cal V}_n \supseteq {\cal V}_{n+1}$.  
By Lemma \ref{test}, $Q$ is a minimal reduction of $A$ with $r_Q(A) = n$ if and only if 
\begin{eqnarray*} \rk M_n(\alpha) & \neq & {n+m-1 \choose m-1},\\
\rk M_{n+1}(\alpha) & = & {n+m \choose m-1}. \end{eqnarray*}
But these two conditions mean $\alpha \in {\cal V}_n \setminus {\cal V}_{n+1}$. \end{pf}

According to Theorem \ref{given},  a number $n \ge 0$ occurs as the reduction number of a minimal reduction of $A$ if and only if ${\cal V}_{n+1}$ is properly contained in ${\cal V}_n$. Since $r(A)$ is the least possible reduction number, $r(A)$ is the least integer $n$ such that ${\cal V}_{n+1} \neq {\cal V}_0$.  By our definition we have $M_0(u) = 0$, so that ${\cal V}_0 = {\Bbb P}_k^{md-1}$. Therefore, $r(A)$ is the largest integer $n$ such that ${\cal V}_n = {\Bbb P}_k^{md-1}$. This fact has the following interesting consequence. \sk

\begin{Corollary} \label{almost}  $r_Q(A) = r(A)$ for almost all minimal reductions $Q$ of $A$.  \end{Corollary}

\begin{pf} As we have seen above, ${\cal V}_{r(A)} = {\Bbb P}_k^{md-1}$. Therefore, $r_Q(A) = r(A)$ if and only if $\alpha \in {\Bbb P}_k^{md-1}\setminus {\cal V}_{r(A)+1}$, which is a non-empty open subset of ${\Bbb P}_k^{md-1}$. \end{pf}

Next we can  characterize the  big reduction number $br(A)$ as follows. Recall that $br(A)$ is defined as the largest possible reduction number of minimal reductions of $A$.\sk
 
\begin{Corollary} \label{big} $br(A)$ is the largest integer $n \ge 0$ such that ${\cal V}_{n+1} \neq {\cal V}_n$.  \end{Corollary}

\begin{pf} This follows from the fact that the possible reduction numbers of minimal reduction reductions of $A$ are exactly the numbers $n$ for which ${\cal V}_{n+1} \neq {\cal V}_n$. \end{pf}

Since $\{{\cal V}_n\}$ is a non-decreasing sequence of varieties, we must have ${\cal V}_{n+1} = {\cal V}_n$ for large $n$. Therefore, from Corollary \ref{big} we obtain a simple proof for the fact that $br(A)$ is finite. In general, $br(A)$ is bounded by the Castelnuovo-Mumford regularity of $A$ [T1, Proposition 3.2]. Bounds for $br(A)$ in terms of the arithmetical degree of $A$ is given in [V3, Proposition 2.2]. It would be of interest to derive these bounds from Theorem \ref{big}.\sk

The big reduction number $br(A)$ can be used to describe the parameter space of Noether normalizations of $A$ which are generated by linear forms.\sk

\begin{Corollary} \label{Noether} Assume that  $A = k[x_1,\ldots,x_m]/I$. Let 
$$y_i = \sum_{j=1}^m\alpha_{ij}x_j,\ i = 1,\ldots,d,$$
where $\alpha = (\alpha_{ij}) \in {\Bbb P}^{md-1}$. Then $k[y_1,\ldots,y_d] \hookrightarrow A$ is a Noether normalization if and only if $\alpha \not\in {\cal V}_{br(A)+1}$. \end{Corollary}

\begin{pf} It is known that $k[y_1,\ldots,y_d] \hookrightarrow A$ is a Noether normalization if and only if the ideal $Q = (y_1,\ldots,y_d,I)/I$ is a minimal reduction of $A$. By Theorem \ref{given}, $Q$ is a minimal reduction of $A$ if and only if there exists an integer $n \ge 1$ such that $\alpha \in {\cal V}_n \setminus  {\cal V}_{n+1}$. Since $\{{\cal V}_n\}$ becomes stationary after $n = br(A)$, this condition means $\alpha \not\in {\cal V}_{br(A)+1}$. \end{pf}

The following example shows that the range of the possible reduction numbers of minimal reductions of $A$ needs not to be a consecutive sequence of integers. \sk

\begin{Example} \label{example1} {\rm Let $2 < a_1 < \cdots < a_{m-1}$ be any sequence of integers. Put $A= k[x_1,\ldots,x_m]/I$, where 
$$I = (x_1^{a_1},\ldots,x_{m-1}^{a_{m-1}})+(x_ix_j|\ 1 \le i < j \le m).$$ Then $\dim A = 1$. Put $a_0 = 2$ and $a_m = \infty$. For $a_{i-1} \le n < a_i$, $i = 1,\ldots,m$, we have $J_n = u_i\cdots u_m$ and hence 
$${\cal V}_n = \{(\alpha_1,\ldots,\alpha_m)|\ \alpha_j = 0\ \text{for some integer}\ i \le j \le m \}. $$
Thus, if we put $Q := (y,I)/I$, where $y = \alpha_1x_1 + \cdots + \alpha_mx_m$, then $Q$ is a minimal reduction of $A$ with
$$r_Q(A) = a_i - 1$$ if $\alpha_{i+1} \neq 0,\cdots,\alpha_m \neq 0$ (cf. [BH, Example 7]). So $1,a_1-1,\ldots,a_{m-1}-1$ is the sequence of the possible reduction numbers of minimal reductions of $A$. In particular, $r(A) = 1$ and $br(A) = a_{m-1}-1$. Moreover, ${\cal V}_{a_{m-1}} = \{(\alpha_1,\ldots,\alpha_m)|\ \alpha_m \neq 0\}$ is the parameter space for $k[y] \hookrightarrow A$ being a Noether normalization. }\end{Example}

Now we shall see that the varieties $\{{\cal V}_n\}$ can be defined without using a presentation of the given graded algebra $A$. This definition will be used later in our study on the core of an ideal. \sk

Assume that the algebra $A$ is generated by the linear forms $a_1,\ldots,a_m$. For every integer $n \ge 0$ we fix a basis of the vector space $A^n$ which consists of monomials of degree $n$ in the elements $a_1,\ldots,a_m$. Let 
$$b_i = \sum_{j=1}^mu_{ij}a_j \  (i = 1,\ldots,d).$$
We write the elements of the form $b_ig$, where $g$ is a monomial in $a_1,\ldots,a_m$ of degree $n-1$, as a linear combination of the elements of the fixed basis of $A_n$ with coefficients in $k[u]$. Let $M'_n(u)$ denote the matrix of these coefficients. Then we denote by $J'_n$ the ideal of $k[u]$ generated by the $h_n\times h_n$ minors of $M'_n(u)$, where $h_n = \dim_kA_n$.  \sk

\begin{Lemma} \label{no-presentation} ${\cal V}_n$ is the projective variety defined by the zero locus of $J'_n$.
\end{Lemma}

\begin{pf} Let $Q = (c_1,\ldots,c_d)$ be an arbitrary ideal generated by $d$ linear forms of $A$. Write 
$$c_i = \sum_{j=1}^m\alpha_{ij}a_j \  (i = 1,\ldots,d),$$
where $\alpha = (\alpha_{ij}) \in {\Bbb P}^{md-1}$. Let $M'_n(\alpha)$ denote the matrix obtained from $M'_n(u)$ by the substitution $u \to \alpha$. It is clear that $Q_n = A_n$ if and only if $\rk M'_n(\alpha) = h_n$. From this it follows that $Q$ is a minimal reduction of $A$ with $r_Q(A) \le n$ if and only if $\alpha$ is not a zero of $I'(n+1)$. Note that $J'_0 = J_0 = 0$. Then using Theorem \ref{given} we can conclude that ${\cal V}_n$ is the zero locus of $J'_n$. \end{pf}

\section{Reduction number and initial ideals}\sk

Let $R=k[x_1,\ldots,x_m]$ be a polynomial ring over a field $k$
and $I$ an arbitrary homogeneous ideal of $R$. Let $\ini(I)$ denote the
initial ideal of $I$ with respect to a term order of $R$. The aim
of this section is to prove that $r(R/I) \le r(R/\ini(I))$. \sk

Given a linear map $\lambda: {\Bbb Z}^m\to {\Bbb Z}$ we can
define a weight order on the monomials of $R$. Let
$\ini_{\lambda}(I)$ denote the initial ideal of $I$ with respect
to this monomial order. It can be shown that $\ini(I) =
\ini_{\lambda}(I)$ for a suitable choice of $\lambda$ (see e.g. [E, p.327]). The ideal
$\ini_{\lambda}(I)$ can be described as follows. \sk

Let $R[t]$ be a polynomial ring over $R$ in one variable $t$. For
any $f \in R[t]$, $f =\sum_ia_iu_i$, where the $u_i$ are monomials
and $0\neq a_i\in k$, we set $b(f)=\max\lambda(u_i)$ and define
\[
f^* :=
t^{b(f)}f(t^{-\lambda(x_1)}x_1,\ldots,t^{-\lambda(x_n)}x_n).
\]
We denote by $I^*$ the ideal of $R[t]$ generated by $\{f^*|\ f\in
I\}$. It is known that $t$ is a non-zerodivisor modulo $I^*$ and that
$R[t]/(I^*,t) \cong R/\ini_{\lambda}(I)$. \sk

In order to study the reduction number of $R/\ini_{\lambda}(I)$ we 
have to pass to the localization $S := R\otimes_k k[t]_{(t)}$ of $R[t]$ which is a standard graded algebra over the local ring $k[t]_{(t)}$ ($\deg t = 0$). Let
$\widetilde{I} = I^*S$. Then $t$ is still a non-zerodivisor modulo $\widetilde{I}$ and
$$S/(\widetilde{I},t) \cong R/\ini_{\lambda}(I).$$

The ideal $\widetilde I$ can be computed from a set of generators of $I$. \sk

\begin{Lemma} \label{generator}
Assume that $I = (f_1,\ldots,f_v)$. Then 
$$\widetilde I = \cup_{n\ge 1}(f_1^*,\ldots,f_v^*)S:t^n.$$
\end{Lemma}

\begin{pf} Since $f_1^*,\ldots,f_v^* \in \widetilde I$ and since $t$ is a non-zerodivisor modulo $\widetilde I$, we have
$$\cup_{n\ge1}(f_1^*,\ldots,f_v^*)S:t^n \subseteq \widetilde I:t^n = \widetilde I.$$
Conversely, let $f$ be an arbitrary element in $I$. Write $f = f_1g_1+\cdots+f_vg_v$ and set $b = \max\{b(f_ig_i)|\ i = 1,\ldots,v\}$. Then
$b(f) \le b$. Hence  
$$t^{b-b(f)}f^* = t^{b-b(f_1g_1)}f_1^*g_1^* + \cdots +
t^{b-b(f_vg_v)}f_v^*g_v^* \in (f_1^*,\ldots,f_v^*)S.$$
From this it follows that $f^* \in (f_1^*,\ldots,f_v^*)S:t^{b-b(f)}$. Since $\widetilde I$ is generated by the elements of the form $f^*$, we get $\widetilde I \subseteq \cup_{n\ge 1}(f_1^*,\ldots,f_v^*)S:t^n$.  \end{pf}\sk

Since $S/\widetilde I$ is a standard graded algebra over a local ring, we can define the reduction number $r(S/\widetilde I)$ of $S/\widetilde I$ as in the case of a standard graded algebra over a field. This reduction number has been used in [HHT] in order to estimate the asymptotic Castelnuovo-Mumford regularity $\reg(\ini_{\lambda}(I^n))$, $n \gg 0$. \sk

\begin{Lemma} \label{tilde} {\rm (i) } A reduction of $S/\widetilde I$ is minimal if and only if it is a reduction of $S/\widetilde I$ generated by $d$ elements, where $d = \dim R/I$, \par
\noindent {\rm (ii) } $r(S/\widetilde I) = r(R/\ini_{\lambda}(I))$ {\rm [HHT, Lemma 3.2]}. \end{Lemma}

\begin{pf} Let $A = S/\widetilde I$ and $\bar A = R/\ini_{\lambda}(I)$. Note that $\bar A = A/tA$ and that $\dim \bar A = d$. Let $t_1,\ldots,t_d$ be arbitrary linear forms in $A$ and $Q = (t_1,\ldots,t_d)$. Let $\bar t_1,\ldots,\bar t_s$
be the images of $t_1,\ldots,t_d$ in $\bar A$ and $\bar Q = (\bar t_1,\ldots,\bar t_d)$. It is clear that if $A_n = Q_n$, then $\bar A_n = \bar Q_n$. Conversely, if $\bar A_n = \bar Q_n$, then $A_n = (Q,t)_n = Q_n + tA_n$. By Nakayama's lemma, this implies $A_n = Q_n$. Thus, $Q$ is a reduction of $A$ if and only if $\bar Q$ is a reduction of $\bar A$ and $r_Q(A) = r_{\bar Q}(\bar A)$ in this case. In particular, $Q$ is a minimal reduction of $A$ if and only if $\bar Q$ is a minimal reduction of $\bar A$. Since $\bar Q$ is a minimal reduction of $\bar A$ if and only if it is a reduction of $\bar A$ which is generated by $d$ elements, $Q$ is a minimal reduction of $A$ if and only if it is a reduction of $A$ which is generated by $d$ elements. Moreover, since $r_Q(A) = r_{\bar Q}(\bar A)$, we can conclude that $r(A) = r(\bar A)$. \end{pf} \sk
 
Using the above observation we are able to prove the following relationships between the reduction numbers of $R/I$ and $R/\ini_{\lambda}(I)$.  \sk

\begin{Theorem} \label{ini-lambda}
For any linear map $\lambda: {\Bbb Z}^m\to {\Bbb Z}$ we have $$r(R/I) \le r(R/\ini_{\lambda}(I)).$$
\end{Theorem}

\begin{pf} Let $T = R \otimes_k k(t)$. By Corollary \ref{extension} 
and Lemma \ref{tilde}(ii) we only need to show that 
$$r(T/IT) \le r(S/\widetilde I).$$
Let $d = \dim R/I$ and $r = r(S/\widetilde I)$. Let $Q$ be a minimal reduction of $S/\widetilde I$ with $r_Q(S/\widetilde I) = r$. By Lemma \ref{tilde}(i) there are $d$ linear forms $y_1,\ldots,y_d$ in $S$ such that $Q = (y_1,\ldots,y_d,\widetilde I)/\widetilde I$. 
We have $$S_{r+1} = (y_1,\ldots,y_d,\widetilde I)_{r+1}.$$\par
Assume that $I = (f_1,\ldots,f_v)$. Note that $T$ is a localization of $S = R \otimes_kk[t]_{(t)}$ and that $t$ is invertible in $T$. Then $\widetilde IT = (f_1^*,\ldots,f_v^*)T$ by Lemma \ref{generator}. Hence 
$$T_{r+1} = (y_1,\ldots,y_d,f_1^*,\ldots,f_v^*)_{r+1}.$$
Let $\phi$ be the graded $k(t)$-automorphism of $T$ with $\phi(x_j) = t^{\lambda(x_j)}x_j$, $j = 1,\ldots,m$. Set $z_i = \phi(y_i)$, $i = 1,\ldots,d$. Since $\phi(f^*) = t^{b(f)}f$ for any $f \in R$, we get
\begin{eqnarray*} T_{r+1}\ =\ \phi(T_{r+1})  
& = & \phi\big((y_1,\ldots,y_d,f_1^*,\ldots,f_v^*)_{r+1}\big)\\ 
& = & (z_1,\ldots,z_d,f_1,\ldots,f_v)_{r+1}\\
& = & (z_1,\ldots,z_d,I)_{r+1}. \end{eqnarray*}
Thus, $(z_1,\ldots,z_d,I)/IT$ is a minimal reduction of $T/IT$ with reduction number not greater than $r$. So we obtain $r(T/IT) \le r = r(S/\widetilde I)$. \end{pf} \sk

The following consequence of Theorem \ref{ini-lambda} gives an affirmative answer to [V2, Conjecture 7.2].  

\begin{Corollary} \label{ini} For any term order of $R$ we have
$$r(R/I) \le r(R/\ini(I)).$$ \end{Corollary}

\begin{pf} We choose an integral weight function $\lambda$ such that
$\ini(I) = \ini_{\lambda}(I)$ and apply Theorem \ref{ini-lambda}.
\end{pf} 

Corollary \ref{ini} has been proved for any generic initial ideal and for Cohen-Macaulay ring $R/I$ by Bresinsky and Hoa [BH, Theorem 12 and Proposition 13]. Moreover, if $\Gin(I)$ is the generic initial ideal of $I$ with respect to the reverse lexicographic order, then $r(R/I) = r(R/\Gin(I))$ by [T2, Theorem 4.3]. \sk

In general, we can not compare $br(R/I)$ and $br(R/\ini(I))$.  First, since $r(R/I)  = r(R/\Gin(I)) = br(R/\Gin(I))$ by [BH, Theorem 11], we have 
$$br(R/I) \ge br(R/\Gin(I)).$$
Since $r(R/I)$ can be strictly less than $br(R/I)$ (see Example \ref{example1}), the above inequality can be a strict inequality. 
On the other hand, if $R/I$ is a Cohen-Macaulay ring, then $br(R/I) = r(R/I)$ by [T, Corollary 3.5] and $r(R/I) \le r(R/\ini(I))$ by Corollary \ref{ini}, hence
$$br(R/I)  \le br(R/\ini(I)).$$
This can be a strict inequality by the following example. 

\begin{Example} \label{example2} {\rm   Let $R = k[x,y,z]$ and $I = (x^2,xz+y^2)$. Then $\ini(I) = (x^2,xz,xy^2,y^4)$ with respect to the lexicographic term order. It is easy to check that
$br(R/I) = 3 < 4 = br(R/\ini(I))$. 
}\end{Example}

The author is grateful to L.~T.~Hoa for showing him the above example.

\section{The local case}\sk

Let $(L,\mm)$ be a local ring with infinite residue field $k$. Let $\aa$ be an ideal of $L$. An ideal $\qq \subseteq \frak a$ is called a reduction of $\aa$ if there is a number $n$ such that $\aa^{n+1} = \qq\aa^n$. In other word, $\qq$ is a reduction of $\aa$ if and only if $\aa$ is integrally dependent on $\qq$. The least number $n$ with the above property will be denoted by $r_\qq(\aa)$. The {\it reduction number} $r(\aa)$ resp. the {\it big reduction number} $br(\aa)$  of $\aa$ is defined as the minimum resp. the maximum of $r_\qq(\aa)$, where $\qq$ is a minimal reduction of $\aa$ with respect to inclusions. \sk

Let $F(\aa)$ denote the {\it fiber ring} $\oplus_{n\ge 0}\aa^n/\mm\aa^n$ of $\aa$. Then $F(\aa)$ is a standard graded algebra over the field $k$. Let $d = \dim F(\aa)$ (the analytic spread of $\aa$). It is known that every minimal reduction of $\aa$ is generated by $d$ elements. Due to [NR], there is the following correspondence between minimal reductions of $\aa$ and minimal reductions of $F(\aa)$. \sk

Let $\qq= (c_1,\ldots,c_d)$ be an ideal generated $d$ elements of $\aa$. Put $Q = (c_1^*,\ldots,c_d^*)$, where $c_1^*,\ldots,c_d^*$ denote the residue classes of $c_1,\ldots,c_d$ in $\aa/\mm\aa$.  Then $\qq$ is a minimal reduction of $\aa$ if and only if $Q$ is a minimal reduction of $F(\aa)$. Moreover, this correspondence preserves the reduction number. \sk
 
\begin{Lemma} \label{fiber} Let $A = F(\aa)$ be the fiber ring of $\aa$. Let $\qq$ be an arbitrary minimal reduction of $\aa$ and $Q$ its corresponding minimal reduction in $A$. Then \par
\noindent {\rm (i)} $r_\qq(\aa) = r_Q(A)$,\par
\noindent {\rm (ii)} $r(\aa) = r(A)$,\par
\noindent {\rm (iii)} $br(\aa) = br(A)$. \end{Lemma}

\begin{pf} We have $A_{n+1} = Q_{n+1}$ if and only if $\aa^{n+1} = \qq\aa^n + \mm\aa^{n+1}$. By Nakayama lemma, the last equation means $\aa^{n+1} = \qq\aa^n$.
Hence
$$r_\qq(A) = \min\{n|\ \aa^{n+1} = \qq\aa^n\} = \min\{n|\  A_{n+1} = Q_{n+1}\} = r_Q(A).$$
By the above correspondence between the minimal reductions of $\aa$ and of $A$, this implies $r(\aa) = r(A)$ and  $br(\aa) = br(A)$.  \end{pf}

Assume that $\aa = (a_1,\ldots,a_m)$.  Let  $L(u)$ denote the local ring $L[u]_{\mm L[u]}$, where $u = \{u_{ij}|\ i = 1,\ldots,d,\ j = 1,\ldots,m\}$ is a family of indeterminates. In $L(u)$ we consider $d$ {\it generic elements} of $\aa$:
$$b_i = \sum_{j=1}^m u_{ij}a_i\ (i = 1,\ldots,d).$$
Put $\aa_u = \aa L(u)$ and $\qq_u = (b_1,\ldots,b_d)L(u)$. 
It is easily seen that $\qq_u$ is a minimal reduction of $\aa_u$. We call $\qq_u$ a {\it generic minimal reduction} of $\aa$.

\begin{Theorem}  Let $\qq_u$ be a generic minimal reduction of $\aa$. Then $$r(\aa) = r_{\qq_u}(\aa_u).$$ \end{Theorem}

\begin{pf} Let $A = F(\aa)$ be the fiber ring of $\aa$. Let $Q_u \subset A_u$ be a generic minimal reduction of  $A$ (in the sense of Section 1). Then $A_u$ is the fiber ring of $\aa_u$ and $Q_u$ is generated by the initial forms of the elements $b_1,\ldots,b_d$ in $A_u$. By Lemma \ref{fiber}(i) we have
$$r_{\qq_u}(\aa_u) = r_{Q_u}(A_u).$$
By Corollary \ref{fiber}(ii)  and Theorem \ref{generic} we have $$r(\aa) = r(A) = r_{Q_u}(A_u).$$
So we obtain $r(\aa) = r_{\qq_u}(\aa_u).$
\end{pf}

Let $\qq= (c_1,\ldots,c_d)$ be an ideal generated $d$ elements of $\aa$. If $$c_i = \sum_{j=1}^m \alpha_{ij}a_j \  (i =1,\ldots,d),$$
we say that $\qq$ is {\it parameterized} by the family $\alpha = (\alpha_{ij}) \in L^{md}$. \sk
 
Let $a^*$ denote the residue class of $a_j$ in $\aa/\mm\aa$ and let $\alpha_{ij}^*$ denote the residue class of $a_{ij}$ in $k = L/\mm$. Then 
$$c_i^* = \sum_{j=1}^m \alpha_{ij}^*a_j ^*\  (i =1,\ldots,d).$$
Hence the corresponding ideal $Q = (c_1^*,\ldots,c_d^*)$ of $\qq$ in $F(\aa)$ is parameterized by the point $\alpha^* = (\alpha_{ij}^*) \in {\Bbb P}_k^{md-1}$.  By passing from $\qq$ to $Q$ we can use the results of Section 2 to describe the parameter space of the minimal reductions of $\aa$ with a given reduction number. In particular, we can compute the big reduction number $br(\aa)$. We leave the readers to formulate the corresponding results. Here we will only describe the parameter space of the minimal reductions of $\aa$. \sk

Fix a minimal basis of the ideal $\aa^{br(\aa)+1}$ which consists of monomials of degree $br(\aa)+1$ in the elements $a_1,\ldots,a_m$. Then we write 
the elements of the form $b_ig$, where $b_i$ is a generic element of $\aa$ ($i = 1,\ldots,d$) and $g$ is a monomial in $a_1,\ldots,a_m$ of degree $br(\aa)$, as a linear combination of the elements of the fixed basis of $\aa^{br(\aa)+1}$ with coefficients in the polynomial ring $L[u]$. Let $M$ denote the matrix of these coefficients. 
Let $\cal J$ denote the ideal of $L[u]$ generated by the $h\times h$-minors of $M$, where $h$ is the minimal number of generators of $\aa^{br(\aa)+1}$. We call $\cal J$ a {\it testing ideal} for the minimal reductions of $\aa$. For every $\alpha \in L^{md}$ let ${\cal J}_\alpha$ denote the ideal of $L$ obtained from $\cal J$ by the specialization $u \to \alpha$. Using this notion we can describe the parameter space of the minimal reductions of $\aa$ as follows. \sk
 
\begin{Proposition} \label{test-local} Let $\qq$ be an ideal generated by $d$ elements of $\aa$ which is parameterized by $\alpha \in L^{md}$. Then $\qq$ is a minimal reduction of $\aa$ if and only if ${\cal J}_\alpha = L$. \end{Proposition}

\begin{pf} Let $A = F(\aa)$ be the fiber ring of $\aa$. Let $Q$ denote the corresponding ideal of $\qq$ in $A$. By \ref{Noether}, $Q$ is a minimal reduction of $A$ if and only if $\alpha^* \not\in {\cal V}_{br(A)+1}$.  Note that $br(A)+1 = br(\aa)+1$ by \ref{fiber}(iii). Let ${\cal J}^*$ denote the ideal  ${\cal J}+\mm L[u]/\mm L[u]$ of the quotient ring $L[u]/\mm L[u] = k[u]$. Then ${\cal J}^*$ is the ideal $J'_n$ introduced at the end of Section 2. 
By Lemma \ref{no-presentation}, ${\cal V}_{br(A)+1}$ is the zero locus of $J'_{br(A)+1} = {\cal J}^*$.  
Therefore, $Q$ is a minimal reduction of $A$ if and only if $\alpha^*$ is not a zero of ${\cal J}^*$. The latter condition means that there exists an element  $f(u) \in {\cal J}$ such that $f(\alpha) \not\in  \mm $. But this is equivalent to say that ${\cal J}_\alpha = L$. \end{pf}

Now we will use the above description of the minimal reductions of $\aa$ to study the core of $\aa$. Recall that the {\it core} of $\aa$, written as $\core(\aa)$, is defined as the intersection of all minimal reductions of $\aa$. As noted by Huneke and Swanson [HS], it is an important question to understand how $\core(\aa)$ is related to $\aa$.  This question is still open. \sk

Since $\aa^{br(\aa)+1} \subseteq \qq$  for any minimal reduction $\qq$ of $\aa$, we have $\aa^{br(\aa)+1} \subseteq \core(\aa).$ 
Between these two ideals we can construct the following ideal. \sk

\begin{Theorem} \label{core} Let $\QQ$ be the ideal of $L[u]$ generated by the generic elements $b_1,\ldots,b_d$ of $\aa$. Let ${\cal J}$ be the testing ideal for the minimal reductions of $\aa$. Then
$$\aa^{br(\aa)+1} \subseteq (\QQ :  {\cal J}^\infty) \cap L \subseteq \core(\aa).$$ \end{Theorem}

\begin{pf} Let $A = F(\aa)$ be the fiber ring of $\aa$. By the definition of the ideal ${\cal J}$ we have ${\cal J}\aa^{br(A)+1} \subseteq \QQ$. Therefore
$$\aa^{br(a)+1} \subseteq (\QQ:{\cal J}) \cap L \subseteq (\QQ :  {\cal J}^\infty) \cap L.$$  
To prove the second inclusion let $x$ be  an arbitrary element of $(\QQ:  {\cal J}^n) \cap L$.  
Then there exists an integer $n$ such that $x{\cal J}^n \subset \QQ$. Let $\qq = (c_1,\ldots,c_d)$ be an arbitrary minimal reduction of $\aa$ which is parameterized by the family $\alpha  \in L^{md}$. By Proposition \ref{test-local} there exits an element $f(u) \in {\cal J}$ such that $f(\alpha)$ is an invertible element in $L$. Since $xf(u)^n \in \QQ$ and since every element of $\QQ$ is specialized to an element of $\qq$ by the substitution $u \to \alpha$, we have $xf(\alpha)^n \in \qq$ so that $x \in \qq$. Therefore, $x \in \core(\aa)$. \end{pf}

Rees and Sally [RS, Theorem 2.6] proved that if $(L,\mm)$ is a Cohen-Macaulay local ring and $\aa$ is an $\mm$-primary ideal, then $\qq_u \cap L \subseteq \core(\aa)$, where  $\qq_u$ is a generic minimal reduction of $\aa$. Recently, Corso, Polini and Ulrich [CPU, Theorem 4.7] showed that this is actually an equality. More generally, they proved that if $L$ is a Cohen-Macaulay local ring and $\aa$ is an universally weakly $(d-1)$-residually $S_2$ ideal satisfying the condition $G_d$, where $d$ is the analytic spread of $\aa$, then $\qq_u \cap L = \core(\aa)$. We refer the readers to their work for the definition of the above notions. Now we shall see that $\qq_u \cap L = (\QQ :  {\cal J}^\infty) \cap L$ in most of the cases considered by Corso, Polini and Ulrich. \sk

\begin{Lemma} \label{contraction}  Let $\qq_u$ be a generic minimal reduction of $\aa$ in $L(u)$. 
Let $\QQ$ and ${\cal J}$ be defined as above.  Assume that every associated prime of $\QQ:(\QQ:\aa^\infty)$ is contained in $\mm L[u]$. Then
$$\qq_u \cap L = (\QQ:  {\cal J}^\infty) \cap L.$$
\end{Lemma}

\begin{pf} Let $\QQ' := \qq_u \cap L[u]$. Since $\qq_u$ is an ideal of the local ring $L(u) = L[u]_{\mm L[u]}$ and since $\qq_u =\QQ L(u)$, every associated prime of $\QQ'$ is contained in $\mm L[u]$. Note that the ideal ${\cal J}+\mm L[u]/\mm L[u]$ of the quotient ring  $L[u]/\mm L[u] = k[u]$ is the ideal $J'_{br(\aa)+1}$ introduced at the end of Section 2. Then, by Theorem \ref{given}(ii) and Lemma \ref{no-presentation} we have  $J'_{br(\aa)+1} \neq 0$. Therefore ${\cal J} \not\in \mm L[u]$. This implies  $\QQ': J = \QQ'$. Hence $\QQ: {\cal J}^\infty \subseteq \QQ'.$ \par
On the other hand, since every associated prime of $\QQ:(\QQ:\aa^\infty)$ is  an associated prime of $\QQ$ contained in $\mm L[u]$, we get $\QQ' \subseteq \QQ:(\QQ:\aa^\infty).$ By the definition of the ideal ${\cal J}$ we have ${\cal J}\aa^{br(A)+1} \subseteq \QQ$. Therefore,  
${\cal J} \subseteq \QQ:\aa^{br(\aa)+1} \subseteq \QQ:\aa^\infty.$ From this it follows that 
$\QQ:(\QQ:\aa^\infty) \subseteq \QQ:{\cal J} \subseteq \QQ:{\cal J}^\infty.$ Hence
$\QQ' \subseteq  \QQ:{\cal J}^\infty.$ \par
So we obtain $\QQ' = \QQ:{\cal J}^\infty$. Hence $\qq_u \cap L  = \QQ' \cap L = (\QQ:{\cal J}^\infty) \cap L$. \end{pf}

It is clear that the afore mentioned result of Rees and Sally is a direct consequence of Theorem \ref{core} and Lemma \ref{contraction}. Combining Lemma \ref{contraction} with the results of [CPU] we obtain the following more general consequence.

\begin{Corollary} Let $L$ be a Cohen-Macaulay local ring and $\aa$  an ideal of $L$ with analytic spread $d$ which satisfies the condition $G_d$. Assume that $\aa$ is universally weakly $(d-1)$-residually $S_2$ and that $\aa L[u]$ is weakly $(d-1)$-residually $S_2$. Then
 $$\core(\aa) = \qq_u \cap L = (\QQ:{\cal J}^\infty) \cap L.$$
\end{Corollary} 

\begin{pf} It was already shown in the proof of [CPU, Proposition 5.4] that under the above assumptions, every associated prime of $\QQ:(\QQ:\aa^\infty)$ is contained in $\mm L[u]$.  Hence the conclusion follows from  Lemma \ref{contraction} and [CPU, Theorem 4.7]. \end{pf}

Corso, Polini and Ulrich [CPU, Example 4.11] also gave an example showing that the formula $\qq_u  \cap L = \core(\aa)$ does not hold for arbitrary ideals in Cohen-Macaulay rings. We shall see that  $(\QQ :  {\cal J}^\infty) \cap L = \core(\aa)$
in their example.

\begin{Example} {\rm Let $L = k[U,V,W]_{(U,V,W)}/(U^2+V^2,VW)$, where $k$ is an infinite field. Denote the images of $U,V,W$ in $R$ by $a_1,a_2,a_3$, respectively. Let $\aa = (a_1,a_2)$. Let $A = F(\aa)$ be the fiber ring of $\aa$. Then  $A \cong  k[U,V]/(U^2+V^2)$. It is easy to check that $br(A) = 1$ so that $br(\aa) = 1$. From this it follows that $\aa^2 \subseteq \core(\aa)$. Since $(a_1)$ and $(a_2)$ are minimal reductions of $\aa$ and since $\aa^2 = (a_1) \cap (a_2)$, we get $\core(\aa) = \aa^2$.  By Theorem \ref{core}, this implies  $(\QQ: {\cal J}^\infty) \cap L = \core(\aa)$. On the other hand, for $L(u) = L[u_1,u_2]_{\mm L[u_1,u_2]}$ and $\qq_u =bL(u)$, we have $\qq_u \cap L = (\aa^2,a_1a_3) \neq \aa^2$. 
}\end{Example}

Due the above observations we raise the question whether $(\QQ :  {\cal J}^\infty) \cap L = \core(\aa)$ holds for arbitrary ideals in arbitrary local rings.

\section*{References}

\noindent [A] I.~M.~Aberbach, Local reduction numbers and Cohen-Macaulayness of associated graded rings, J. Algebra 178 (1995), 833-842.\par

\noindent [AHT] I.~M.~Aberbach, C.~Huneke and N.~V.~Trung, Reduction numbers, Brian\c con-Skoda theorems and depth of the Rees algebras, Compositio Math. 97 (1995), 403-434. \par

\noindent [BH] H.~Bresinsky and L.~T.~Hoa, On the reduction number of some graded algebras, Proc. Amer. Math. Soc. 127 (1999), 1257-1263.\par

\noindent [C] A. Conca, Reduction numbers and initial ideals, preprint.\par

\noindent [CPU] A.~Corso, C.~Polini and B. Ulrich, The structure of the core of ideals, Math. Ann. (to appear). \par

\noindent [E] D. Eisenbud, Commutative algebra, with a view
toward algebraic geometry, Springer, 1994. \par

\noindent [ES] P. Eakin and A. Sathaye, Prestable ideals, J. Algebra 41 (1976), 439-454. \par

\noindent [GNN] S.~Goto, Y.~Nakamura and K.~Nishida, Cohen-Macaulay graded rings associated to ideals, Amer. J. Math. 118 (1996), 1196-1213.\par

\noindent [JU] M.~Johnson and B. Ulrich, Artin-Nagata properties and Cohen-Macaulay associated graded rings, Compositio Math. 108 (1996), 7-29.

\noindent [JK] B. Johnston and D. Katz, Castelnuovo regularity and graded rings associated to an ideal, Proc. Amer. Math. Soc. 123 (1995), 727-734.\par

\noindent [HHT] J. Herzog, L.T. Hoa and N.V. Trung, Asymptotic linear bounds for the Castelnuovo-Mumford regularity, Preprint, 2000. \par 

\noindent [Hu] S.~Huckaba, Reduction numbers for ideals of higher analytic spread, Math. Proc. Camb. Phil. Soc. 102 (1987), 49-57. \par

\noindent [HH1] S.~Huckaba and C.~Huneke, Powers of ideals having small analytic deviation, Amer. J. Math. 114 (1992), 367-403.\par

\noindent [HH2] S.~Huckaba and C.~Huneke, Rees algebras of ideals having small analytic deviation, Trans. Amer. Math. Soc. 339 (1993), 373-402.\par

\noindent [HS] C. Huneke and I. Swanson, Cores of ideals in 2-dimensional regular local rings, Michigan Math. J. 42 (1995), 193-208. \par

\noindent [L] J. Lipman, Equimultiplicity, reductions, and blowing up, in: Commutative algebra, analytic method, 111--147, Lect. Notes in Pure and Appl. Math. 68, Dekker, 1979.\par

\noindent [NR] D. G. Northcott and D. Rees, Reductions of ideals
in local rings, Proc. Cambridge Philos. Soc. 50 (1954), 145-158.\par

\noindent [RS] D. Rees and J. Sally, General elements and joint reductions,
Michigan Math. J. 35 (1988), 241-254. \par

\noindent [S] J. Sally,  Tangent cones at Gorenstein singularities, Compositio Math. 40 (1980), 167-175.\par

\noindent [T1] N.~V.~Trung, Reduction exponent and degree bound
for the defining equations of graded rings, Proc. Amer. Math.
Soc. 101 (1987), 229-236. \par

\noindent [T2] N.~V.~Trung, Gr\"obner bases, local cohomology and
reduction number, Proc. Amer. Math. Soc. 129 (1) (2001), 9-18.\par

\noindent [U] B.~Ulrich, Artin-Nagata properties and reductions of ideals, Contemporary Math. 159 (1994), 373-400. \par

\noindent [V1] W. Vasconcelos, The reduction number of an algebra,
Compositio Math. 106 (1996), 189-197. \par

\noindent [V2] W. Vasconcelos, Cohomological degrees of graded
modules, in: Six lectures on commutative algebra (Bellaterra, 1996),
345--392, Progr. Math. 166, Birkh\"auser, Basel, 1998. \par

\noindent [V3] W. Vasconcelos, Reduction numbers of ideals, 
J. Algebra 216 (1999), 652-664. \par

\end{document}